\documentclass[11pt]{article}
\usepackage{amsmath}
\usepackage{amssymb}
\usepackage{amsthm}
\usepackage{mathtools}
\usepackage{amssymb}
\usepackage{mathrsfs}
\usepackage[all]{xy}
\usepackage{mathrsfs}
\usepackage{indentfirst}
\usepackage{color}
\usepackage[unicode=true,colorlinks=true,linkcolor=blue,citecolor=blue,urlcolor=blue]{hyperref}
\allowdisplaybreaks
\newtheorem{theorem}{Theorem}[section]
\newtheorem*{theorem*}{Theorem}
\newtheorem{corollary}[theorem]{Corollary}
\newtheorem{remark}[theorem]{Remark}
\newtheorem{lemma}[theorem]{Lemma}
\newtheorem{definition}[theorem]{Definition}

\newtheorem{proposition}[theorem]{Proposition}
\newtheorem{example}[theorem]{Example}
\newtheorem*{proposition*}{Proposition}

\newcommand{\R}{\mathbb{R}}

\newcommand{\be}{\begin{eqnarray*}}
\newcommand{\ee}{\end{eqnarray*}}
\newcommand{\ba}{\begin{align*}}
\newcommand{\bpm}{\begin{pmatrix}}
\newcommand{\epm}{\end{pmatrix}}

\begin{document}
\title{On the Fueter--Sce--Qian theorem for partial-slice monogenic functions}
\author{Qinghai Huo$^1$\thanks{This work was partially supported by Fundamental Research Funds for the Central Universities (No. JZ2025HGTB0171) and China Scholarship Council (No. 202506690055).}
,\ Irene Sabadini$^2$\thanks{This work was partially supported by PRIN 2022 {\em Real and Complex Manifolds: Geometry and Holomorphic Dynamics. I.S. is a member of GNSAGA of INdAM}.}\ ,
 and Zhenghua Xu$^1$\thanks{This work was partially supported by  Fundamental Research Funds for the Central Universities (No. JZ2025HGTG0250) and China Scholarship Council (No. 202506690052).}\\
 \\
\emph{$^1$\small School of Mathematics, Hefei University of Technology,}
\emph{\small  Hefei, 230601, P.R. China}\\
\emph{\small E-mail address:  2022800145@hfut.edu.cn;  zhxu@hfut.edu.cn}
\\
\emph{\small $^2$Dipartimento di Matematica, Politecnico di Milano,}
\emph{\small Via E. Bonardi, 9, 20133 Milano, Italy} \\
\emph{\small E-mail address:   irene.sabadini@polimi.it}
}

\maketitle

\begin{abstract}
The notion of partial-slice monogenic functions  has recently been introduced to provide a unified framework encompassing both monogenic and slice monogenic functions over Clifford algebras. In this new setting, we establish some  fractional Laplacian identities for the Cauchy kernels and their derivatives in the sense of tempered distributions. These identities allow us to  prove a termwise Fueter--Sce--Qian  theorem   for partial-slice monogenic functions represented by convergent Laurent series consisting only of principal parts.
\end{abstract}
\section{Introduction}
 In 1934,  R. Fueter   provided a method to produce  regular  functions of one quaternionic variable from  holomorphic  functions of one complex variable  \cite{Fueter}.
 In the   language of slice analysis,  the Laplacian in $\mathbb{R}^{4}$
maps the slice regular function  induced by a holomorphic  function into an axially regular  function, see \cite{Colombo-Sabadini-Sommen-10}.
In 1957, the Fueter  theorem  was generalized by Sce \cite{Sce}   to higher dimensions  $\mathbb{R}^{n+1}$ with odd $n\in \mathbb{N}$. In modern terms, for odd $n\in \mathbb{N}$,  the power $\Delta^{\frac{n-1}{2}}$ of  the Laplacian   $\Delta $   in $\mathbb{R}^{n+1}$ maps the slice monogenic function  induced by a holomorphic  function into an axially monogenic  function, see \cite{Colombo-Sabadini-Sommen-10}.
 By using the Fourier multiplier definition of the fractional Laplacian,  Qian  generalized   the results of Fueter and Sce   to the case of even $n$, see \cite{Qian-97} and \cite{Qian-01}. In particular, Qian defined the monomial functions in $\mathbb{R}^{n+1}$
\begin{equation*}
P^{(-k)}( x)=\tau [ x^{-k} ], \quad P^{(k-1)}(x)=I[ P^{(-k)}](x),  \quad k\in \mathbb{N}\setminus \{0\},\end{equation*}
 where  $\tau=\Delta^{\frac{n-1}{2}}$  is defined  by  the Fourier transformation $\mathcal{F}$ in   $\mathbb{R}^{n+1}$ as
$$\tau f(x)=\mathcal{F}^{-1} \big[  (2\pi i| \xi|)^{n-1} \mathcal{F}[f](\xi) \big](x),$$
and  $I$ is the (left) Kelvin inversion   given by
$$I[f]( x)=\frac{\overline{ x}}{ | x|^{n+1} }f( x^{-1}).$$
It should be pointed out that  $x^{-k}$ in  the definition of $P^{(-k)}$ should be understood as a tempered distribution for $k\geq n+1$, see Example \ref{x^k} and Definition \ref{T} below.

Nowadays this deep result,  known as the Fueter--Sce--Qian theorem, has become one of the fundamental results in hypercomplex analysis.
The main purpose of this work is to establish fractional Laplacian identities for Cauchy-type kernels and their derivatives in the sense of tempered distributions, and to derive a termwise Fueter--Sce--Qian theorem for certain  \textit{partial-slice monogenic functions}, originally called  \textit{generalized partial-slice monogenic functions}. This class of functions was  introduced  over Clifford algebras in 2023 \cite{Xu-Sabadini} (see also \cite{Xu-Sabadini-3}) to  include the two theories of monogenic functions  and of slice monogenic functions as special cases. In this new setting, the theory is rapidly developing and several results  have been obtained, such as the conformal invariance \cite{Ding-Xu}, the  Fueter--Sce theorem \cite{Xu-Sabadini-2,Xie}, an Almansi-type decomposition \cite{Huo-24}, a Segal--Bargmann transform \cite{Xu-Sabadini-4}, a Teodorescu transform \cite{Hu-Ding}.   Partial-slice monogenic functions   have  been  further investigated in the non-associative algebra  of octonions  \cite{Xu-Sabadini-26} and in the more general setting of real alternative $\ast$-algebras \cite{Huo-Sabadini-Xu-26}.

  The concept of partial-sliceness has also been developed by Ghiloni and Stoppato by introducing the theory of $T$-regular functions in \cite{Ghiloni-Stoppato-24-1,Ghiloni-Stoppato-24-2}, where  a unified theory of regularity in one hypercomplex variable for an alternative $\ast$-algebra was proposed, including slice regular functions \cite{Gentili-Struppa-07} and partial-slice monogenic functions as special cases.  In particular,  the Fueter--Sce phenomenon for $T$-regular functions over general associative $\ast$-algebras was investigated  \cite{Ghiloni-Stoppato-25}.  Besides, Perotti showed  that   a suitable power of the Laplacian maps Dunkl-regular functions to Dunkl monogenic functions with axial symmetries, which gives the so-called Fueter trees. Some of these Fueter trees are   described by partial-slice monogenic functions in  the context of octonions  and  Clifford algebras \cite{Perotti-26}. For more versions of the Fueter theorem    in various settings, see e.g.  \cite{Altavilla,De-Adan,Dong-Qian-21,Eelbode,Fei,Kou,Perotti-22,Qian-15,Sommen-00}.

Very recently, Colombo, De Martino and Debernardi Pinos \cite{Colombo-Martino} have developed  a method for defining (distributionally) the  fractional Laplacian $(-\Delta)^{\frac{1}{2}}$ for smooth functions with power-type singularities at the origin, including the negative monomials $x^{-k-1}$ for $x\in\mathbb{R}^{n+1}\setminus\{0\}$ and $k\in \mathbb{N}$, where the test function space is  a subspace  of the Lizorkin space (i.e., the intersection of the Lizorkin space and its image under the Fourier transform, denoted by $\Xi(\mathbb{R}^{n+1})$). In particular,  their method gives   the  fractional Laplacian $(-\Delta)^{\frac{1}{2}}$ on a class of holomorphic functions whose Laurent series only has principal part  and derives  several applications, with particular emphasis on  the quaternionic fine structures. The fine structures, see \cite{CDP25,Fivedim}, are obtained by refining the Fueter map, namely  the Laplacian (in the quaternionic case) or its powers (in the Clifford algebra case), and then factorizing it by means of the Dirac operator and its conjugate and associating the corresponding spaces of functions in their kernels.  In this construction, instead of considering the Dirac operator and its conjugate one make take fractional powers of the Laplacian, and this is the approached followed in \cite{Colombo-Martino}.

In the present  paper, we study the fractional Laplacian $(-\Delta)^s$ on the  Cauchy-type  kernels
 \[ K_\lambda(x)=\frac{\overline{x}}{|x|^{\lambda+2}}, \qquad x\in\mathbb{R}^{n+1}\setminus\{0\}, -1<\lambda<n,\]
and their  derivatives of arbitrary order,  for all possible  $s\in \mathbb{R}$ satisfying   $-1<\lambda+2s<n$  in the sense of tempered distributions  on the Schwartz  class, instead of the test space  $\Xi(\mathbb{R}^{n+1})$ used in \cite{Colombo-Martino}.
The distributional Fourier transforms of the kernels are locally integrable functions, and the products with the fractional Laplacian multipliers remain locally integrable under the stated parameter restrictions. Thus, the fractional Laplacians of these kernels are well defined in the space of tempered distributions via the Fourier transform. Based on these distributional identities,   we can  prove a  termwise Fueter--Sce--Qian theorem for some partial-slice monogenic functions,   as special cases,  generalizing  the corresponding results in \cite{Colombo-Martino} for quaternionic slice regular functions   and in \cite{Colombo-Martino-Qian-S} for slice monogenic  functions.

The plan of the article is as follows. Section 2 collects preliminary notions on Clifford algebras and partial-slice monogenic functions, including Fueter polynomials and Cauchy kernels. In particular, Theorem \ref{thm:partial-slice-kernel} provides further examples of partial-slice monogenic functions. Section 3 establishes fractional Laplacian identities for Cauchy-type kernels and their derivatives  in the sense of tempered distributions.  Based on these  distributional identities,  in Section 4 we prove  a termwise Fueter--Sce--Qian theorem for partial-slice monogenic functions of type $(p,q)$   represented by convergent Laurent series consisting only of principal parts.  To this end, we establish  a  sharp mixed derivative estimate in Proposition \ref{prop:sharp-mixed-derivative-estimate} by   the Banach theorem.

\section{Preliminaries}
In this section, we collect some preliminary results on Clifford algebras and on   partial-slice monogenic
functions recently introduced in \cite{Xu-Sabadini}.
\subsection{Clifford algebras}
As usual, denote by $\mathbb{N}$ all  nonnegative  integers. For $n\in  \mathbb{N}\setminus \{0\},$ let $\{e_1,e_2,\ldots, e_n\}$ be the standard orthonormal basis for the   real Euclidean space  $\mathbb{R}^n$. A real Clifford algebra, denoted by $\mathbb{R}_{n}$, is generated  by these basis elements assuming that   $$e_i e_j + e_j e_i= -2\delta_{ij}, \quad 1\leq i,j\leq n,$$
where $\delta_{ij}$ is the Kronecker symbol.
Each element in  $\mathbb{R}_{n}$ can be written as
 $$a=\sum_A a_Ae_A, \quad a_A\in \mathbb{R},$$
 where $e_A=e_{j_1}e_{j_2}\cdots e_{j_r},$
$A=\{j_1,j_2, \cdots, j_r\}\subseteq\{1, 2, \cdots, n\}$ with $1\leq j_1< j_2 < \cdots < j_r \leq n,$ and  set $ e_\emptyset=e_0=1$.

As a real  vector space, the dimension of the Clifford algebra $\mathbb{R}_n$  is $2^{n}$. The norm of $a$ is defined by $|a|= ({\sum_{A}|a_{A}|^{2}} )^{\frac{1}{2}}.$ For $k=0,1,\ldots,n$, the real linear subspace $\mathbb{R}_n^k$ of $\mathbb{R}_n$, called $k$-vectors, is  generated by the $\begin{pmatrix} n\\k\end{pmatrix}$ elements of the form
 $$e_A=e_{i_1}e_{i_2}\cdots e_{i_k},\quad 1\leq i_1<i_2<\cdots<i_k\leq n.$$
In particular, the element $a_{0}=a_{\emptyset}$ is called the scalar part of $a$.

The Clifford conjugation of $a$ is  defined as
$$\overline{a} =\sum_Aa_A\overline{e_A},$$
 where $\overline{e_{j_1}\ldots e_{j_r}}=\overline{e_{j_r}}\ldots\overline{e_{j_1}},\ \overline{e_j}=-e_j,1\leq j\leq n,$ and $ \overline{e_0}=e_0=1$.

The set of  paravectors of $\mathbb{R}_n$ is   $\mathbb{R}_n^{0} \oplus \mathbb{R}_n^{1}$, which shall  be identified with $\mathbb{R}^{n+1}$ via
$$(x_0,x_1,\ldots,x_n) \longmapsto    x=\sum_{i=0}^{n}e_ix_i.$$
For a paravector $ x\neq0$, its norm is given by $| x|=( x\overline  x)^{1/2}$  and so its inverse is given by $ x^{-1}=  \overline{ x}| x|^{-2}.$
\begin{proposition}(\cite[Theorem 3.14]{Gurlebeck}) \label{paravectors}
 For all  $a\in\mathbb R_n$  and paravectors $x\in\mathbb R^{n+1}$, it holds that
\[
|xa|=|ax|=|x|\,|a|.\]
\end{proposition}
\subsection{Partial-slice monogenic functions}
Let $n\in  \mathbb{N}\setminus \{0\}$ and   consider functions $f:\Omega\longrightarrow  \mathbb{R}_{n}$, where $\Omega $ is an open set in $\R^{n+1}$. From now on, let $p\in \{0,1,\ldots,n-1\}, q=n-p$, and    split the paravector  $x$ in the Clifford algebra $\mathbb{R}_{n}$ as
$$x=\boldsymbol{x}_p+\underline{\boldsymbol{x}}_q \in\R^{p+1}\oplus\R^q, \quad \boldsymbol{x}_p=\sum_{i=0}^{p}x_i e_i,\ \underline{\boldsymbol{x}}_q=\sum_{i=p+1}^{p+q}x_i e_i.$$
Similarly, the generalized Cauchy--Riemann operator $D$ is split as
$$D=D_{\boldsymbol{x}_p} +D_{\underline{\boldsymbol{x}}_q}, \quad D_{\boldsymbol{x}_p} =\sum_{i=0}^{p}e_i\partial_{x_i}, \
D_{\underline{\boldsymbol{x}}_q} =\sum_{i=p+1}^{p+q}e_i\partial_{x_i}.$$

Denote by $\mathbb{S}$ the sphere of unit $1$-vectors in $\mathbb R^q$, whose elements $(x_{p+1},\ldots, x_{p+q})$ are identified with $\underline{\boldsymbol{x}}_q=\sum_{i=p+1}^{p+q}x_i e_i$, i.e.
$$\mathbb{S}=\big\{\underline{\boldsymbol{x}}_q: \underline{\boldsymbol{x}}_q^2 =-1\big\}=\big\{\underline{\boldsymbol{x}}_q=\sum_{i=p+1}^{p+q}x_i e_i:\sum_{i=p+1}^{p+q}x_i^{2}=1\big\}.$$

Note that, for $\underline{\boldsymbol{x}}_q\neq0$, there exist  uniquely  $r\in \mathbb{R}^{+}=\{x\in \mathbb{R}: x>0\}$ and $\underline{\omega}\in \mathbb{S}$, such that $\underline{\boldsymbol{x}}_q=r\underline{\omega}$, more precisely
 $$r=|\underline{\boldsymbol{x}}_q|, \quad \underline{\omega}=\frac{\underline{\boldsymbol{x}}_q}{|\underline{\boldsymbol{x}}_q|}. $$
 When $\underline{\boldsymbol{x}}_q= 0$, set $r=0$ and $\underline{\omega}$ is not uniquely defined, since $x=\boldsymbol{x}_p+\underline{\omega} \cdot 0$ for every $\underline{\omega}\in \mathbb{S}$.

The closed upper half-space $\mathrm{H}_{\underline{\omega}}$ in $\mathbb{R}^{p+2}$ associated with $\underline{\omega}\in \mathbb{S}$ is defined by
$$\mathrm{H}_{\underline{\omega}}=\{\boldsymbol{x}_p+r\underline{\omega}, \boldsymbol{x}_p \in\R^{p+1}, r\geq0 \},$$
and it is clear from the previous discussion that
$$ \R^{p+q+1}=\bigcup_{\underline{\omega}\in \mathbb{S}} \mathrm{H}_{\underline{\omega}}, \quad \R^{p+1}=\bigcap_{\underline{\omega}\in \mathbb{S}} \mathrm{H}_{\underline{\omega}}.$$

Now we recall  the  concept of  partial-slice monogenic functions.
\begin{definition}[Partial-slice monogenic functions] \label{definition-slice-monogenic}
 Let $\Omega$ be an open set in $\mathbb{R}^{n+1}$. A function $f :\Omega \rightarrow \mathbb{R}_{n}$ is called left  partial-slice monogenic of type $(p,q)$ if, for all $ \underline{\omega} \in \mathbb S$, its restriction $f_{\underline{\omega}}$ to
 $\Omega_{\underline{\omega}}:=\Omega\cap (\mathbb{R}^{p+1} \oplus \underline{\omega} \mathbb{R})\subseteq \mathbb{R}^{p+2}
$  has continuous partial derivatives and  satisfies
$$D_{(p,q)}f(x):=(D_{\boldsymbol{x}_p}+\underline{\omega}\partial_{r}) f_{\underline{\omega}}(\boldsymbol{x}_p+r\underline{\omega})=0,$$
for all $x=\boldsymbol{x}_p+r\underline{\omega} \in \Omega_{\underline{\omega}}$.
 \end{definition}
From now on, we denote by $\mathcal{SM}_{(p,q)}(\Omega)$ (or $\mathcal{SM}_{(p,q)}^{L}(\Omega)$ when needed) the  right Clifford module of  all  left partial-slice monogenic functions  of type $(p,q)$  on $\Omega$. This class of functions  was  originally referred to as \textit{generalized partial-slice monogenic functions}  of type $(p,q)$ on $\Omega$ and  denoted by $\mathcal{GSM}(\Omega)$ in \cite{Xu-Sabadini,Xu-Sabadini-3}.

In particular,  $\mathcal{SM}_{(n-1,1)}(\Omega)$ reduces  to the class of  left   monogenic functions on $\Omega$,   denoted by $\mathcal{M}(\Omega)$. Meanwhile, $\mathcal{SM}_{(0,n)}(\Omega)$ reduces  to the class of  left  slice monogenic functions on $\Omega$,   denoted by $\mathcal{SM}(\Omega)$.
See \cite{Brackx} and \cite{Colombo-Sabadini-Struppa-09}.

Likewise, denote by $\mathcal{SM}_{(p,q)}^{R}(\Omega)$ the left Clifford module of  all right  partial-slice monogenic functions of type $(p,q)$ $f:\Omega  \rightarrow \mathbb{R}_{n}$ which are defined by requiring that the restriction $f_{\underline{\omega}}$ to $\Omega_{\underline{\omega}}$ satisfies
$$f(x)D_{(p,q)}:={f_{\underline{\omega}} (\boldsymbol{x}_p+r\underline{\omega})D_{\boldsymbol{x}_p}}+ \partial_{r}f_{\underline{\omega}} (\boldsymbol{x}_p+r\underline{\omega})\underline{\omega}=0, \quad x=\boldsymbol{x}_p+r\underline{\omega} \in \Omega_{\underline{\omega}},$$
for all  $ \underline{\omega} \in \mathbb S$.

In the setting of partial-slice monogenic functions,  Fueter variables and Fueter polynomials are defined as follows.
\begin{definition}\label{Fueter-variables}
 The   (left) Fueter variables   are defined as
 $$ z_{\ell}(x)=z_{\ell}^{L}(x)= z_{\ell}(\boldsymbol{x}_p+r\underline{\omega})= x_{\ell}+r \underline{\omega} e_\ell, \ \ell=0,1,\ldots,p.$$
Similarly, the  right  Fueter variables   are defined as
 $$ z_{\ell}^{R}(x)= z_{\ell}^{R}(\boldsymbol{x}_p+r\underline{\omega})=  x_{\ell}+r e_\ell  \underline{\omega},\ \ell=0,1,\ldots,p.$$
 \end{definition}

For $\mu=(\mu_0,\mu_1,\ldots,\mu_p)\in\mathbb{N}^{p+1}$, denote
$$|\mu|=\mu_0+\mu_1+\cdots+\mu_p,\quad \mu!=\mu_0!\ldots\mu_p!.$$
\begin{definition}[Fueter polynomials]\label{definition-Fueter}
For $\mu=(\mu_0,\ldots,\mu_p)\in\mathbb N^{p+1}$, let $\overrightarrow{\mu}=(j_1,\ldots,j_k)$ be the nondecreasing sequence in which each $\ell\in\{0,\ldots,p\}$ occurs exactly $\mu_\ell$ times, where $k=|\mu|$. Define
$$\mathcal{P}^{L}_{\mu}(x)=\frac{1}{k!}\sum_{(i_1,i_2,\ldots,i_k)\in\sigma(\overrightarrow{\mu})}z_{i_1}z_{i_2}\cdots z_{i_k},$$
where the sum runs over the $\dfrac{k!}{\mu!}$ different permutations $\sigma(\overrightarrow{\mu})$ of $\overrightarrow{\mu}$. When $\mu=(0,\ldots,0)=\mathbf{0}$, we set $\mathcal{P}^{L}_{\mathbf{0}}(x)=1$; when there is at least one negative component in $\mu$, we set $\mathcal{P}^{L}_{\mu}(x)=0$.

Similarly, we can define $\mathcal{P}^{R}_{\mu}$ when $z_{\ell}$ are replaced by $z_{\ell}^{R}$.
\end{definition}

\begin{definition} \label{slice-domain}
 Let $\Omega$ be a domain in $\mathbb{R}^{p+q+1}$.

(i)   $\Omega$ is called  slice domain if $\Omega\cap\mathbb R^{p+1}\neq\emptyset$  and $\Omega_{\underline{\omega}}$ is a domain in $\mathbb{R}^{p+2}$ for every  $\underline{\omega}\in \mathbb{S}$.

(ii)  $\Omega$   is called  partially  symmetric with respect to $\mathbb R^{p+1}$  if, for   $\boldsymbol{x}_p\in\R^{p+1}, r \in \mathbb R^{+},$ and $ \underline{\omega}  \in \mathbb S$,
$$x=\boldsymbol{x}_p+r\underline{\omega} \in \Omega\Longrightarrow [x]:=\boldsymbol{x}_p+r \mathbb S=\{\boldsymbol{x}_p+r \underline{\omega}, \ \  \underline{\omega}\in \mathbb S\} \subseteq \Omega. $$
 \end{definition}

  Partial-slice monogenic  functions  have a  Representation formula \cite[Theorem~3.14]{Xu-Sabadini}.
\begin{theorem}  {\bf(Representation Formula)}  \label{Representation-Formula-SM}
Let $\Omega\subseteq \mathbb{R}^{p+q+1}$ be a  partially symmetric slice domain and $f:\Omega\rightarrow \mathbb{R}_{p+q}$ be a  partial-slice monogenic function.  Then, for any $\underline{\omega}\in \mathbb{S}$ and for $\boldsymbol{x}_p+r\underline{\omega} \in \Omega$,
$$
f(\boldsymbol{x}_p+r \underline{\omega})=\frac{1}{2} (f(\boldsymbol{x}_p+r\underline{\eta} )+f(\boldsymbol{x}_p-r\underline{\eta}) )+
\frac{ 1}{2} \underline{\omega}\underline{\eta} (  f(\boldsymbol{x}_p-r\underline{\eta} )-f(\boldsymbol{x}_p+r\underline{\eta})),
$$
for any $\underline{\eta}\in \mathbb{S}$.

Moreover, the following two functions do not depend on $\underline{\eta}$:
$$F_1(\boldsymbol{x}_p,r)=\frac{1}{2} (f(\boldsymbol{x}_p+r\underline{\eta} )+f(\boldsymbol{x}_p-r\underline{\eta} ) ),$$
$$F_2(\boldsymbol{x}_p,r)=\frac{ 1}{2}\underline{\eta}(  f(\boldsymbol{x}_p-r\underline{\eta} )-f(\boldsymbol{x}_p+r\underline{\eta})).$$
\end{theorem}

We introduce the kernel function \cite[Example 3.6]{Xu-Sabadini} with a brief calculation.
\begin{example}\label{Cauchy-kernel-example}
Let $\Omega=\mathbb{R}^{n+1}\setminus \{0\}$ and denote
$$E(x) := \frac{\overline{x}}{|x|^{p+2}}, \quad x \in\Omega.$$
Then $E\in \mathcal{SM}_{(p,q)}^{L}(\Omega) \cap{\mathcal{SM}}_{(p,q)}^{R}(\Omega).$
\end{example}
For $x=\boldsymbol{x}_p+r\underline\omega\ne0$, we have $D_{(p,q)}\overline{x}=p+2$ and $x\overline{x}=|x|^2$. Hence
\[
D_{(p,q)}E(x)=(p+2)|x|^{-p-2}-(p+2)|x|^{-p-4}x\overline{x}=0.
\]
The right case follows analogously.

For $y\in \mathbb{R}^{n+1}$, denote
\[
E_y(x):=E(y-x)=\frac{\overline{y-x}}{|y-x|^{p+2}}, \quad x\in \mathbb{R}^{n+1}\backslash \{y\}.
\]
Generally,  $E_y$ is not  partial-slice monogenic of type $(p,q)$.
Hence, we introduce the Cauchy kernel \cite[Definition 3.31]{Xu-Sabadini}  as follows.
\begin{definition} Given $y \in \mathbb{R}^{n+1}$,  define  the \textit{left partial-slice monogenic Cauchy kernel}
\begin{equation*}
 \mathcal{E}_{y}^{L}(x)=\frac{1}{2}(1-\underline{\omega}\underline{\eta})E_{y}(\pi_{y}(x))+
\frac{1}{2}(1+\underline{\omega}\underline{\eta})E_{y}( \pi_{y}(x)_{\diamond}),
\end{equation*}
where $x=\boldsymbol{x}_p+r  \underline{\omega},$ $y=\boldsymbol{y}_p+\widetilde{r}  \underline{\eta},
\pi_{y}(x)=\boldsymbol{x}_p+r \underline{\eta} $  and $\pi_{y}(x)_{\diamond}=\boldsymbol{x}_p-r\underline{\eta}$.

Likewise,   define the \textit{right partial-slice monogenic  Cauchy kernel}
 \begin{equation*} \mathcal{E}_{y}^{R}(x)=\frac{1}{2}E_{y}(\pi_{y}(x))(1-\underline{\eta}\underline{\omega})+
\frac{1}{2}E_{y}( \pi_{y}(x)_{\diamond})(1+\underline{\eta}\underline{\omega}).
\end{equation*}
\end{definition}

The   Cauchy kernel $\mathcal E^{L}_{y}(x)$ is well-defined on $\mathbb R^{n+1}\setminus[y]$. Indeed, if $r=0$, the defining formula reduces to $\mathcal E^{L}_{y}(\boldsymbol{x}_p)=E_{y}(\boldsymbol{x}_p)$, independently of the choice of $\underline{\omega}$. Similarly, if $\widetilde r=0$, a direct computation gives
\[
\mathcal E^{L}_{\boldsymbol{y}_p}(x)=E_{\boldsymbol{y}_p}(x)= \frac{\overline{\boldsymbol{y}_p-x}}{|\boldsymbol{y}_p-x|^{p+2}},
\]
which is independent of the choice of $\underline{\eta}$. Consequently, the definition is unambiguous even when either $x$ or $y$ belongs to $\mathbb R^{p+1}$. Furthermore, we have
$$\mathcal{E}^{L}_{y}\in \mathcal{SM}_{(p,q)}^{L}(\mathbb{R}^{n+1}\setminus [y]), \ \mathcal{E}_{y}^{R}\in \mathcal{SM}_{(p,q)}^{R}(\mathbb{R}^{n+1}\setminus [y]).$$
When $(p,q)=(0,n)$, $\mathcal{E}^{L}_{y}(x)$ reduces to the left slice monogenic Cauchy kernel
$$S_{L}^{-1}(y,x)
= \big(x^{2}-2\,\mathrm{Re}(y)x+|y|^{2}\big)^{-1}\big(\overline{y }-x\big).$$
Meanwhile,  $\mathcal{E}_{y}^{R}(x)$ reduces to the right slice monogenic Cauchy kernel
\[
S_{R}^{-1}(y,x)=  (\overline{y}-x)\big(x^{2}-2\,\mathrm{Re}(y)x+|y|^{2}\big)^{-1}.
\]

Let $\lambda\in\mathbb{R}$, $y\in\mathbb{R}^{n+1}$ and define
$$\Omega_y=\bigl\{x\in\mathbb{R}^{n+1}:y^2-2x_0y+|x|^2\in\mathbb{R}^{n+1}\setminus(-\infty,0]\bigr\}.$$
Consider the function
$$k_{\lambda}(y,x):=(y-\bar{x})(y^2-2x_0y+|x|^2)^{-\lambda},\qquad x\in\Omega_y,$$
We clarify the meaning of the power $Q^{-\lambda}$  with $Q=y^2-2x_0y+|x|^2$. Write $y=y_0+\underline y$. If $\underline y\neq0$, set $I=\underline y/|\underline y|$ and $\mathbb C_I=\mathbb R+I\mathbb R$. Then
\[
Q=y^2-2x_0y+|x|^2=\bigl(y_0^2-|\underline y|^2-2x_0y_0+|x|^2\bigr)+2(y_0-x_0)\underline y\in\mathbb C_I,
\]
and define $Q^{-\lambda}:=\exp(-\lambda\operatorname{Log}Q)$ using the principal branch $\operatorname{Log}$ of the logarithm in defined on $\mathbb{C}_I\setminus (-\infty,0]$.  If $y\in\mathbb R$, then $Q^{-\lambda}=|x-y|^{-2\lambda}$ is the usual negative real power for $x\neq y$.

\begin{theorem}\label{thm:partial-slice-kernel}
For $\lambda\in\mathbb{R}$ and $y\in\mathbb{R}^{n+1}$,  we have
\begin{equation*}
D_{(p,q)} k_{\lambda}(y,x)=(2\lambda-p-2)(y^2-2x_0y+|x|^2)^{-\lambda}.
\end{equation*}
In particular,    $k_{\lambda}(y, \cdot)$ is left partial-slice monogenic of type $(p,q)$   if and only if
$\lambda=\frac{p}{2}+1.$
\end{theorem}
\begin{proof}
Write $x=x_0+\underline{x}$, where $\underline{x}=\sum_{j=1}^{p}x_je_j+r\underline{\omega}$, and set
\[
P=y-\bar{x},\qquad Q=y^2-2x_0y+|x|^2.
\]
Then $k_{\lambda}(y,x)=PQ^{-\lambda}$. Since $Q$ and its partial derivatives commute, the usual differentiation rule for powers gives
\[
\partial_{x_0}Q^{-\lambda}=-2\lambda(x_0-y)Q^{-\lambda-1},\qquad
\partial_{x_j}Q^{-\lambda}=-2\lambda x_jQ^{-\lambda-1},
\]
for $j=1,\ldots,p$, and
\[
\partial_rQ^{-\lambda}=-2\lambda rQ^{-\lambda-1}.
\]
Moreover, using $e_j^2=\underline{\omega}^2=-1$, we obtain
\[
D_{(p,q)}P=-(p+2).
\]
Consequently, applying the product rule yields
\begin{align*}
D_{(p,q)}k_{\lambda}(y,x)
=-(p+2)Q^{-\lambda} -2\lambda\bigl[P(x_0-y)+\underline{x}P\bigr]Q^{-\lambda-1}.
\end{align*}
On the other hand, since $P=-(x_0-y)+\underline{x}$ and $\underline{x}^2=-|\underline{x}|^2$, we have
\begin{align*}
P(x_0-y)+\underline{x}P
&=\bigl(-(x_0-y)+\underline{x}\bigr)(x_0-y)
+\underline{x}\bigl(-(x_0-y)+\underline{x}\bigr)\\
&=-(x_0-y)^2+\underline{x}^2\\
&=-\bigl(y^2-2x_0y+|x|^2\bigr)=-Q.
\end{align*}
It follows that
\[D_{(p,q)}k_{\lambda}(y,x)=(2\lambda-p-2)Q^{-\lambda}.
\]
Hence, we conclude that
\[
D_{(p,q)}k_{\lambda}(y,x)=0\quad\Longleftrightarrow\quad \lambda=\frac{p}{2}+1,
\]
which  completes the proof.
\end{proof}
As a special case, we have  \cite[Theorem 2.15]{Colombo-Martino-Qian-S}:
\begin{corollary}The function $k_{\lambda}(y,x)$ is left monogenic in the variable $x$ if and only if $\lambda=\frac{n+1}{2}$.
\end{corollary}

 \section{Fueter--Sce--Qian  theorem for kernel functions}
In this section, we realize the Cauchy-type kernels and their derivatives as tempered distributions and establish their fractional Laplacian identities by the Fourier transformation.

 Denote  by  $\mathscr{S}(\mathbb{R}^{n+1})$ the Schwartz space of rapidly decreasing  functions on $\mathbb{R}^{n+1}$.
\begin{definition}\label{Fourier-transform}
For $f\in \mathscr{S}(\mathbb{R}^{n+1})$, define
the Fourier transform and inverse Fourier transform of   $f$ as
$$\mathcal{F}[f]( \xi) = \widehat{f}( \xi)=    \int _{\mathbb{R}^{n+1}} e^{2\pi i \langle x, \xi\rangle} f(x) dx,$$
and
$$\mathcal{F}^{-1}[f](x) = f^\vee (x)=  \int _{\mathbb{R}^{n+1}} e^{-2\pi i \langle x, \xi\rangle} f(\xi) d\xi, $$
respectively.

For  a tempered distribution $T\in \mathscr{S}'(\mathbb{R}^{n+1})$,  define
its Fourier transform and inverse Fourier transform of   $T$ by
\[
\langle \widehat{T},\varphi\rangle=\langle T,\widehat{\varphi}\rangle,\
\langle T^\vee,\varphi\rangle=\langle T,\varphi^\vee\rangle, \quad \varphi\in \mathscr{S}(\mathbb{R}^{n+1}).
\]
\end{definition}
Throughout this paper, the test function spaces $\mathscr S(\mathbb R^{n+1})$ and $\Phi(\mathbb R^{n+1})$ (see Definition \ref{Lizorkin} below) consist of complex-valued functions. Let $\mathbb C_n=\mathbb R_n\otimes_{\mathbb R}\mathbb C$, where the imaginary unit $i$ of $\mathbb{C}$ commutes with every element of $\mathbb R_n$. Fourier transforms and distributional operations on Clifford-valued functions and distributions are understood componentwise. In particular, a $\mathbb C_n$-valued tempered distribution has the form $T=\sum_A T_Ae_A$, with $T_A\in\mathscr S'(\mathbb R^{n+1})$, and its Fourier transform is defined by $\widehat T=\sum_A\widehat{T_A}e_A$. We use the same notation $\mathscr S'(\mathbb R^{n+1})$ and $\Phi'(\mathbb R^{n+1})$ for these distribution spaces.

\begin{remark}
For $s\in\mathbb R$, put $m_s(\xi)=(2\pi|\xi|)^{2s}$.  By the Fourier multiplier, the fractional Laplacian
$$(-\Delta)^sT=\mathcal F^{-1}[m_s\widehat T],\qquad T\in\mathscr S'(\mathbb R^{n+1}),$$
requires the product $m_s\widehat T$ to be well defined as a tempered distribution. For $s\in\mathbb N$, the multiplier $m_s$ is a polynomial, so the definition applies to every tempered distribution. For noninteger $s>0$, the multiplier is not smooth at the origin; for $s<0$, it is singular there. Thus, for $s\notin\mathbb N$, the product must be justified for the distribution under consideration.
\end{remark}

One standard framework uses the Lizorkin space of test functions; see \cite{Rubin2025,Samko}.
\begin{definition}[Fractional Laplacian on Lizorkin spaces]\label{Lizorkin}
Let
$$\Phi(\mathbb R^{n+1})=\left\{\varphi\in\mathscr S(\mathbb R^{n+1}):\int_{\mathbb R^{n+1}}x^\alpha\varphi(x)\,dx=0\text{ for every }\alpha\in\mathbb N^{n+1}\right\}.$$
Equip $\Phi(\mathbb R^{n+1})$ with the topology induced by $\mathscr S(\mathbb R^{n+1})$, and denote its continuous dual by $\Phi'(\mathbb R^{n+1})$.
For $s\in\mathbb R$ and $\varphi\in\Phi(\mathbb R^{n+1})$, define
$$(-\Delta)^s\varphi=\mathcal F^{-1}\!\left[(2\pi|\xi|)^{2s}\widehat\varphi(\xi)\right].$$
The infinite-order vanishing of $\widehat\varphi$ at the origin ensures
$$(-\Delta)^s (\Phi(\mathbb R^{n+1}))=\Phi(\mathbb R^{n+1}).$$
 For $T\in\Phi'(\mathbb R^{n+1})$, define by duality
$$\langle(-\Delta)^sT,\varphi\rangle=\langle T,(-\Delta)^s\varphi\rangle,\qquad \varphi\in\Phi(\mathbb R^{n+1}).$$
\end{definition}

\begin{remark}\label{Fueter-laplace}
Every $\mathbb R_n$-valued polynomial $P$ represents the zero distribution in $\Phi'(\mathbb R^{n+1})$ by the vanishing moment condition. Consequently,
$$(-\Delta)^sP=0\qquad\text{in }\Phi'(\mathbb R^{n+1}),\qquad s\in\mathbb R.$$
In particular, this applies to the Fueter polynomials $\mathcal P_\mu^L$ and $\mathcal P_\mu^R$, $\mu\in\mathbb N^{p+1}$. Thus the Lizorkin distribution space does not distinguish these polynomials.

In $\mathscr S'(\mathbb R^{n+1})$, however, nonzero polynomials represent nonzero distributions. Let $(p,q)=(0,n)$, where $n\ge2$ is even, and set $s=(n-1)/2$. The Fourier transform of $x^k$, $k\in\mathbb N$, is a linear combination of derivatives of $\delta_0$ of order $k$. For $0\le k\le n-2$, the multiplier $(2\pi|\xi|)^{n-1}$ is sufficiently differentiable to multiply these distributions, and all its derivatives of order at most $k$ vanish at the origin. Hence
$$(-\Delta)^{(n-1)/2}x^k=0\qquad\text{in }\mathscr S'(\mathbb R^{n+1}),\qquad 0\le k\le n-2.$$
For $k\ge n-1$, this multiplication is not justified by the usual multiplication rule for derivatives of $\delta_0$, since $|\xi|^{n-1}$ is not of class $C^{n-1}$ at the origin. An additional prescription is therefore needed in these cases.
\end{remark}

We now turn to negative powers and Cauchy-type kernels, for which we specify tempered-distribution extensions.

\begin{example}\label{x^k}
For $k,n\in\mathbb N\setminus\{0\}$, consider the slice monogenic function
\[
f_k(x):=x^{-k},
\qquad
x\in\mathbb{R}^{n+1}\setminus\{0\}.
\]
Since
$|x^{-k}|=|x|^{-k},$
we have, in polar coordinates,
$$\int_{|x|<1}|x^{-k}|\,dx
=|\mathbb S^n|\int_0^1 r^{n-k}\,dr.$$
where $\mathbb S^n=\{u\in\mathbb R^{n+1}:|u|=1\}$ and $|\mathbb S^n|=2\pi^{(n+1)/2}/\Gamma((n+1)/2)$ denotes its surface area.
Then
\[
x^{-k}\in L^1_{\mathrm{loc}}(\mathbb{R}^{n+1})
\quad\Longleftrightarrow\quad
k<n+1.
\]
Hence,  for
$k=1,\ldots,n,$
the function \(f_k\) induces a   tempered distribution by
\[
\langle f_k,\varphi\rangle=\int_{\mathbb{R}^{n+1}}x^{-k}\varphi(x)\,dx,
\qquad \varphi\in  \mathscr{S}(\mathbb{R}^{n+1}).
\]
For \(k\geq n+1\),  \(x^{-k}\) does not induce a  regular  tempered distribution by integration.
\end{example}
Nevertheless, this does not mean that \(x^{-k}\) cannot be
realized as a tempered distribution. In fact, it admits tempered distribution extensions, for instance by the Hadamard finite-part construction or by analytic continuation of homogeneous distributions (see e.g. \cite[p. 43-45]{Samko}). In the present paper, we use the canonical extension as follows.
\begin{definition}\label{T}
Let $T_{-1}$ be   the (regular) tempered distribution induced by the
locally integrable function \(x^{-1}=\overline{x}/|x|^2\), and define
\[
T_{-k}:=\frac{(-1)^{k-1}}{(k-1)!}\partial_{x_0}^{\,k-1}T_{-1}\in\mathscr{S}'(\mathbb{R}^{n+1}),\qquad k\geq1.
\]
\end{definition}
Note  that $T_{-k}$ is homogeneous of degree $-k$, and its restriction to $\mathbb{R}^{n+1}\setminus\{0\}$ is represented by the smooth function $x^{-k}$.

More generally, for $\lambda\in \mathbb{R}$, denote
$$K_\lambda(x):=\frac{\overline{x}}{|x|^{\lambda+2}}, \qquad x\in\mathbb{R}^{n+1}\setminus\{0\}.$$
\begin{lemma}\label{Kerel-distribution}Let $n\in\mathbb N\setminus\{0\}$. Then,
for $\lambda<n$,  $K_\lambda$  defines a  regular  tempered distribution on $\mathbb{R}^{n+1}$.
\end{lemma}

\begin{proof}
Since $|\overline{x}|=|x|,$
we have
\[|K_\lambda(x)|=\frac{|\overline{x}|}{|x|^{\lambda+2}}=|x|^{-\lambda-1}.
\]
In a neighborhood of the origin we obtain
$$\int_{|x|<1}|K_\lambda(x)|\,dx=|\mathbb S^n|\int_0^1 r^{n-\lambda-1}\,dr.$$
Hence, we have for  $\lambda<n$
\[
\int_{|x|<1}|K_\lambda(x)|\,dx<\infty,
\]
and  then
\[
K_\lambda\in L^1_{\mathrm{loc}}(\mathbb{R}^{n+1}).
\]
For $|x|\geq1$, we have
\[
|K_\lambda(x)|=|x|^{-\lambda-1} \leq(1+|x|)^{\max\{0,-\lambda-1\}}.
\]
Thus $K_\lambda$ has at most polynomial growth at infinity and defines a regular tempered distribution on $\mathbb R^{n+1}$.
\end{proof}

For any multi-index
$
\mu=(\mu_0,\ldots,\mu_n)
\in\mathbb{N}^{n+1},$
we write
\[|\mu|=\mu_0+\cdots+\mu_n,
\qquad
\mu!=\mu_0!\cdots\mu_n!.
\]

\begin{definition}\label{T-distribution}
Let  $\lambda<n$ and  $\mu \in\mathbb{N}^{n+1}$.  Consider the  tempered distribution  $K_\lambda \in \mathscr{S}'(\mathbb{R}^{n+1})$, and define
\[
T_{\mu}^{(\lambda)}
:=
\frac{(-1)^{|\mu|}}{\mu!}\,
\partial^\mu K_\lambda \in \mathscr{S}'(\mathbb{R}^{n+1}).
\]
  In particular, denote
\[
T_{-k}^{(\lambda)}:=T_{(k-1,0,\ldots,0)}^{(\lambda)},
\qquad T_{-k}:=T_{-k}^{(0)},\qquad k\geq1.
\]
\end{definition}

 To prove main  theorems, we need a technical lemma \cite[p. 73]{Stein}.
\begin{lemma} \label{Stein}
Let  $n\in \mathbb{N}$  and $0<\alpha<n+1$.
Suppose that $H_k(x)$ is a homogeneous harmonic polynomial of degree $k\geq1$. Then it holds that for   $\varphi \in \mathscr{S}(\mathbb{R}^{n+1})$
$$\int_{\mathbb{R}^{n+1}} \frac{H_k(x)}{|x|^{k+n+1-\alpha}}\widehat{\varphi}(x)dx=\gamma_{k,\alpha} \int_{\mathbb{R}^{n+1}} \frac{H_k(x)}{|x|^{k+\alpha}} \varphi(x)dx,$$
where
\begin{equation}
\label{constant}\gamma_{k,\alpha}:=i^{k}\pi^{\frac{n+1}{2} -\alpha}  \frac{\Gamma((k+\alpha)/2)}{\Gamma((k+n+1-\alpha)/2)}.
\end{equation}
\end{lemma}

\begin{remark}  In $\mathscr{S}'(\mathbb{R}^{n+1})$,  Lemma \ref{Stein} implies that
\begin{equation*}
\mathcal{F}\Big[\frac{H_k(x)}{|x|^{k+n+1-\alpha}}\Big](\xi) =\gamma_{k,\alpha}  \frac{H_k(\xi)}{|\xi|^{k+\alpha}}, \end{equation*}
and
\begin{equation*}
\frac{H_k(x)}{|x|^{k+n+1-\alpha}} =\gamma_{k,\alpha}  \mathcal{F}^{-1}\Big[\frac{H_k(\xi)}{|\xi|^{k+\alpha}}\Big](x).\end{equation*}
 Although the function $H_k(x)$   is real (or complex),  it can be  directly generalized to  the     Clifford algebras valued case.
 \end{remark}

\begin{lemma} \label{Kp-Fourier}
Let  $-1< \lambda<n$,  $\mu \in\mathbb{N}^{n+1}$.   Then we have   in $\mathscr{S}'(\mathbb{R}^{n+1})$
\begin{equation} \label{T-mu-p}
\widehat{T_\mu^{(\lambda)}}(\xi)
=\frac{\gamma_{1,n-\lambda}}{\mu!} (2\pi i\xi)^\mu
\frac{\overline{\xi}}{|\xi|^{n+1-\lambda}},
\end{equation}
where $\gamma_{1,n-\lambda}$ is given by (\ref{constant}).
\end{lemma}
\begin{proof} For $-1< \lambda<n$, we apply Lemma \ref{Stein} with
\[
k=1, \qquad
H_1(x)=\overline{x}, \qquad \alpha=n-\lambda\in (0,n+1).
\]
Hence
\[
\frac{H_1(x)}
{|x|^{k+n+1-\alpha}}
=
\frac{\overline{x}}{|x|^{\lambda+2}}
=
K_\lambda(x).
\]
Therefore, in the sense of tempered distributions,
\begin{equation} \label{K-Fourier}
\widehat{K_\lambda}(\xi)=\gamma_{1,n-\lambda}\frac{\overline{\xi}}{|\xi|^{n+1-\lambda}}.\end{equation}

Since
\[\widehat{\partial^\mu K_\lambda} =(-2\pi i\xi)^\mu\widehat{K_\lambda},
\]
we obtain
$$\widehat{T_\mu^{(\lambda)}}
=
\frac{(-1)^{|\mu|}}{\mu!}
(-2\pi i\xi)^\mu
\widehat{K_\lambda}=
\frac{(2\pi i\xi)^\mu}{\mu!}
\widehat{K_\lambda}.$$
Substituting (\ref{K-Fourier}), we get
$$\widehat{T_\mu^{(\lambda)}}(\xi)
=\frac{\gamma_{1,n-\lambda}}{\mu!}
(2\pi i\xi)^\mu
\frac{\overline{\xi}}
{|\xi|^{n+1-\lambda}},$$
which completes the proof. \end{proof}

For $s\in\mathbb R$, denote $m_s(\xi):=(2\pi|\xi|)^{2s}$. In view of Lemma~\ref{Kp-Fourier}, we define
\begin{equation}\label{definition-laplace}
\Big\langle m_s\widehat{T_\mu^{(\lambda)}},\varphi\Big\rangle
:=\frac{\gamma_{1,n-\lambda}}{\mu!}
\int_{\mathbb R^{n+1}}(2\pi i\xi)^\mu\overline{\xi}|\xi|^{\lambda-n-1}m_s(\xi)\varphi(\xi)\,d\xi,
\quad \varphi\in\mathscr S(\mathbb R^{n+1}).
\end{equation}

\begin{proposition}\label{well-defined}
Let $-1<\lambda<n$, $s\in\mathbb R$, and $\mu\in\mathbb N^{n+1}$. If $\lambda+2s>-1$, then $m_s\widehat{T_\mu^{(\lambda)}}$ defines a tempered distribution through \eqref{definition-laplace}.
\end{proposition}

\begin{proof}Let $-1<\lambda<n$, $s\in\mathbb R$, and $\mu\in\mathbb N^{n+1}$.  We have
\begin{equation}\label{laplace-right}
m_s(\xi)\widehat{T_\mu^{(\lambda)}}(\xi)
=\frac{\gamma_{1,n-\lambda}}{\mu!}(2\pi)^{2s}
(2\pi i\xi)^\mu\overline{\xi}
|\xi|^{\lambda+2s-n-1},
\end{equation}
where the equality holds pointwise for $\xi\ne0$ through the function representation in Lemma~\ref{Kp-Fourier}.

Since
$$\big|m_s(\xi)\widehat{T_\mu^{(\lambda)}}(\xi)\big|
\le C_{\lambda,\mu,s}|\xi|^{|\mu|+\lambda+2s-n},$$
the assumption $\lambda+2s>-1$ gives
$$\int_0^1r^{|\mu|+\lambda+2s}\,dr<\infty.$$
Thus, the function in \eqref{laplace-right} is locally integrable and has at most polynomial growth at infinity. Consequently, the integral in \eqref{definition-laplace} defines a tempered distribution, i.e.
$$m_s\widehat{T_\mu^{(\lambda)}}\in\mathscr S'(\mathbb R^{n+1}).$$
The proof is complete.
\end{proof}

Since $\mathcal F^{-1}$ is a continuous linear automorphism of $\mathscr S'(\mathbb R^{n+1})$, the fractional Laplacian $(-\Delta)^sT_\mu^{(\lambda)}$ is well-defined by
$$(-\Delta)^sT_\mu^{(\lambda)}
:=\mathcal F^{-1}\Big(m_s\widehat{T_\mu^{(\lambda)}}\Big)
\in\mathscr S'(\mathbb R^{n+1}).$$

\begin{theorem}\label{p-p+2s}
Let $-1< \lambda<n$ and $s\in \mathbb{R}$ satisfy
$-1<\lambda+2s<n.$
For every multi-index
$\mu\in\mathbb N^{n+1}$, we have
\[
(-\Delta)^s T_\mu^{(\lambda)}=c_{n,\lambda,s}\,T_\mu^{(\lambda+2s)},
\quad\text{in }\mathscr S'(\mathbb R^{n+1}),
\]
 where
\begin{equation} \label{constant-distribution}
c_{n,\lambda,s}:=(2\pi)^{2s}\frac{\gamma_{1,n-\lambda}}{\gamma_{1,n-\lambda-2s}}= 2^{2s}
\frac{\Gamma\left(\frac{n-\lambda+1}{2}\right)
\Gamma\left(\frac{\lambda+2s+2}{2}\right)}{\Gamma\left(\frac{\lambda+2}{2}\right)
\Gamma\left(\frac{n-\lambda-2s+1}{2}\right)}.\end{equation}
\end{theorem}
\begin{proof}Let $-1< \lambda<n$.
By (\ref{laplace-right}),  we have
$$(2\pi|\xi|)^{2s}
\widehat{T_\mu^{(\lambda)}}(\xi)
=\frac{\gamma_{1,n-\lambda}}{\mu!}(2\pi)^{2s}
(2\pi i\xi)^\mu
\frac{\overline{\xi}}
{|\xi|^{n+1-\lambda-2s}}.$$
On the other hand,  formula (\ref{T-mu-p})
applied with $\lambda$ replaced by $\lambda+2s\in (-1,n)$ gives
\[
\widehat{T_\mu^{(\lambda+2s)}}(\xi)
=
\frac{\gamma_{1,n-\lambda-2s}}{\mu!}
(2\pi i\xi)^\mu
\frac{\overline{\xi}}
{|\xi|^{n+1-\lambda-2s}}.
\]
Hence
\[
(2\pi|\xi|)^{2s} \widehat{T_\mu^{(\lambda)}}(\xi)=c_{n,\lambda,s} \widehat{T_\mu^{(\lambda+2s)}}(\xi),
\]
where
\[
c_{n,\lambda,s}
=
(2\pi)^{2s}\frac{\gamma_{1,n-\lambda}} {\gamma_{1,n-\lambda-2s}}.
\]
By   definition, we have
\[
\mathcal F[(-\Delta)^sT_\mu^{(\lambda)}]=(2\pi|\xi|)^{2s}\widehat{T_\mu^{(\lambda)}}.
\]
Since the Fourier transform is a continuous linear automorphism of
$\mathscr{S}'(\mathbb R^{n+1})$, we conclude that
\[
(-\Delta)^sT_\mu^{(\lambda)}=c_{n,\lambda,s}T_\mu^{(\lambda+2s)},
\]
which  proves the result.
\end{proof}

By Theorem \ref{p-p+2s} for $s=\pm \frac{1}{2}$ and $\mu=0$,  we have:
\begin{corollary} \label{Kernel-relation}
For $p=0,1,2,\ldots,n-2$, we have
$$(-\Delta)^{\frac{1}{2}} K_p =2\pi \frac{\gamma_{1,n-p}}{\gamma_{1,n-p-1}} K_{p+1},
\quad\text{in }\mathscr S'(\mathbb R^{n+1}).$$
\end{corollary}
\begin{corollary} \label{Kernel-relation-inverse}
For $p=1,2,\ldots,n-1$, we have
$$(-\Delta)^{-\frac{1}{2}} K_{p} =\frac{1}{2\pi} \frac{\gamma_{1,n-p}}{\gamma_{1,n-p+1}}K_{p-1},
\quad\text{in }\mathscr S'(\mathbb R^{n+1}).$$
\end{corollary}
In Corollaries \ref{Kernel-relation} and  \ref{Kernel-relation-inverse}, upon identifying the left-hand side with its smooth representative away from the origin, the identity holds pointwise on $\mathbb R^{n+1}\setminus\{0\}$.

As a direct consequence of Theorem \ref{p-p+2s}, we can recover   \cite[Theorem 5.1]{Colombo-Martino-Qian-S}.
Set $c_{0,n}:=c_{n,0,(n-1)/2}=2^{n-1}\Gamma\left(\frac{n+1}{2}\right)^2$.

\begin{corollary}\label{cor:translated-Fueter-kernel}
Let $n\geq2$ and $y\in\mathbb R$. Then we have in $\mathscr S'(\mathbb R^{n+1})$
$$(-\Delta)^{\frac{n-1}{2}}S_L^{-1}(y,x)=c_{0,n}(y-\overline{x})(y^2-2yx_0+|x|^2)^{-(n+1)/2}.$$
\end{corollary}
\begin{proof}
Taking $\lambda=0$, $\mu=0$, and $s=(n-1)/2$ in Theorem \ref{p-p+2s}, we obtain
\[
(-\Delta)^{\frac{n-1}{2}}K_0=c_{0,n}K_{n-1},
\qquad c_{0,n}=2^{n-1}\Gamma\Big(\frac{n+1}{2}\Big)^2.
\]
For $y\in\mathbb R$, identified  with $(y,0,\ldots,0)\in\mathbb R^{n+1}$, we have $$ |x-y|^2=y^2-2yx_0+|x|^2,$$ and
\[S_L^{-1}(y,x)=-(\tau_yK_0)(x), \]
where  $\tau_y$ denotes the translation  $(\tau_y K_0)(x)=K_0(x-y)$.
By the Fourier transform convention in Definition~\ref{Fourier-transform}, we have
\[
\widehat{\tau_yT}(\xi)=e^{2\pi i y\xi_0}\widehat T(\xi),\qquad T\in\mathscr S'(\mathbb R^{n+1}),\]
where $\langle\tau_yT,\varphi\rangle=\langle T,\varphi(\,\cdot+y)\rangle$.
Hence
\[
(2\pi|\xi|)^{n-1}\widehat{\tau_yK_0}(\xi)=e^{2\pi i y\xi_0}(2\pi|\xi|)^{n-1}\widehat K_0(\xi)
=c_{0,n}e^{2\pi i y\xi_0}\widehat K_{n-1}(\xi),
\]
which    defines a tempered distribution.
Define
\[(-\Delta)^{\frac{n-1}{2}}\tau_yK_0:=\mathcal{F}^{-1}[m_{(n-1)/2}\widehat{\tau_yK_0}].\]
Taking the inverse Fourier transform  yields
\[
(-\Delta)^{\frac{n-1}{2}}\tau_yK_0
=\tau_y(-\Delta)^{\frac{n-1}{2}}K_0
=c_{0,n}\tau_yK_{n-1}, \quad {\rm in} \  \mathscr S'(\mathbb R^{n+1}),
\]
i.e.
\[
(-\Delta)^{\frac{n-1}{2}}S_L^{-1}(y,x)
=-c_{0,n}K_{n-1}(x-y)=c_{0,n}\frac{y-\bar x}{|x-y|^{n+1}},
\]
 as asserted.
\end{proof}

\section{Fueter--Sce--Qian theorem for some partial-slice monogenic functions}
Ignoring the nonzero multiplicative constants, the preceding result in Corollary  \ref{Kernel-relation} gives the {following chain of (fractional) Fueter--Sce--Qian maps:}
$$
K_0\xmapsto{\,(-\Delta)^{1/2}\,}K_1
\xmapsto{\,(-\Delta)^{1/2}\,}K_2
\xmapsto{\,(-\Delta)^{1/2}\,}\cdots
\xmapsto{\,(-\Delta)^{1/2}\,}K_{n-1}.
$$
{Recalling  Example \ref{Cauchy-kernel-example}, the chain of maps and kernels obtained above motivates  to conjecture the existence of the chain}
\[
 \mathcal{SM}
\xrightarrow{\;(-\Delta)^{\frac12}\;}
\mathcal{SM}_{(1,n-1)}
\xrightarrow{\;(-\Delta)^{\frac12}\;}
\mathcal{SM}_{(2,n-2)}
\cdots
\xrightarrow{\;(-\Delta)^{\frac12}\;}
\mathcal{SM}_{(n-2,2)}
\xrightarrow{\;(-\Delta)^{\frac12}\;}
\mathcal{M}.
\]
More precisely, we conjecture that
\[
(-\Delta)^{\frac12}:\mathcal{SM}_{(p,q)}\longrightarrow\mathcal{SM}_{(p+1,q-1)},\qquad p=0,1,\ldots, n-2,
\]
under suitable assumptions on the functions and on the distributional realization of the fractional Laplacian.

For parameters $\beta,\gamma>-1$, the Jacobi polynomial
$P_m^{(\beta,\gamma)}$ is defined as
\[
P_m^{(\beta,\gamma)}(y)=\frac{(\beta+1)_m}{m!}\sum_{k=0}^{m}\binom{m}{k}\frac{(m+\beta+\gamma+1)_k}{(\beta+1)_k}\big(\frac{y-1}{2}\big)^k,\qquad y\in[-1,1],
\]
where  $(a)_0=1$ and $(a)_k=a(a+1)\cdots(a+k-1)$ for $k\geq1$.

In \cite[Proposition 3.2]{Colombo-Martino},   a formula for the partial derivatives of the function $|x|^{-\alpha}$ was obtained in terms of Jacobi polynomials, which can be used to  estimate   the derivatives \cite[Proposition 3.4]{Colombo-Martino}.

\begin{proposition}\label{c-d}
Let $x\in\mathbb{R}^{n+1}\setminus\{0\}$, $\alpha>0$, and $m\in\mathbb{N}\setminus \{0\}$. Then, we have
\[\frac{\partial^{2m}}{\partial x_j^{2m}} \frac{1}{|x|^\alpha}
=\frac{(-1)^m\sqrt{\pi}\,
\Gamma\big(m+\frac{\alpha}{2}\big)(2m)!} {|x|^{2m+\alpha}
\Gamma\big(\frac{\alpha}{2}\big) \Gamma\big(m+\frac12\big)} P_m^{\big(-\frac12,\frac{\alpha-1}{2}\big)}
\Big( 1-\frac{2x_j^2}{|x|^2}\Big),\]
and
\[\frac{\partial^{2m+1}}{\partial x_j^{2m+1}}\frac{1}{|x|^\alpha}=\frac{(-1)^{m+1}\sqrt{\pi}\,
\Gamma\big(m+\frac{\alpha}{2}+1\big)(2m+1)!x_j}{|x|^{2m+2+\alpha}\Gamma\big(\frac{\alpha}{2}\big)
\Gamma\big(m+\frac32\big)}P_m^{\big(\frac12,\frac{\alpha-1}{2}\big)}
\Big(1-\frac{2x_j^2}{|x|^2}\Big),\]
for all $j=0,1,\ldots,n$.  In particular, for any $\ell\in\mathbb N\setminus \{0\}$ and $\alpha\geq2$, we have
\[
\left|\frac{\partial^\ell}{\partial x_j^\ell}\frac{1}{|x|^\alpha}\right|\leq C_\alpha\frac{\ell^{\alpha-1}\ell!}{|x|^{\ell+\alpha}}.
\]
\end{proposition}
By the Banach theorem, we can establish an optimal  mixed derivative estimate without using an explicit Jacobi polynomial formula.
\begin{proposition}\label{prop:sharp-mixed-derivative-estimate}
Let $\alpha>0$ and $\beta\in\mathbb N^{n+1}$, and put $\ell=|\beta|$. Then
\[
\left|\partial^\beta |x|^{-\alpha}\right|
\leq(\alpha)_\ell |x|^{-\alpha-\ell},\qquad x\in\mathbb R^{n+1}\setminus\{0\},
\]
where $(\alpha)_\ell=\Gamma(\alpha+\ell)/\Gamma(\alpha)$   is optimal uniformly over all multi-indices of order $\ell$ and all nonzero $x$.

In particular, there exists $C_\alpha>0$, independent of $n$, $\beta$, and $x$, such that
\[\left|\partial^\beta |x|^{-\alpha}\right|
\leq C_\alpha(1+\ell)^{\alpha-1}\ell! |x|^{-\alpha-\ell}, \qquad x\in\mathbb R^{n+1}\setminus\{0\}.
\]\end{proposition}
\begin{proof}
The case $\ell=0$ is immediate. Assume $\ell\geq1$ and set $f(x)=|x|^{-\alpha}$. Fix $x\in\mathbb R^{n+1}\setminus\{0\}$ and a unit vector $v\in\mathbb R^{n+1}$, and put
$$R=|x|,\qquad a=\frac{\langle x,v\rangle}{R}\in[-1,1].$$
By the generating function for Gegenbauer polynomials \cite[p.~82]{Szego}
$$(1-2au+u^2)^{-\lambda} =\sum_{k=0}^{\infty}C_k^\lambda(a)u^k, \qquad |u|<1, \lambda>0,$$
we have, for real $t$ with $|t|<R$,
$$f(x-tv)=R^{-\alpha}\big(1-2a\frac{t}{R}+\frac{t^2}{R^2}\big)^{-\alpha/2}
=\sum_{k=0}^{\infty}C_k^{\alpha/2}(a)R^{-\alpha-k}t^k.$$
Consequently,
$$\mathrm d^\ell f(x)[v,\ldots,v]=(-1)^\ell\ell!R^{-\alpha-\ell}C_\ell^{\alpha/2}(a), $$
where $\mathrm d^\ell f(x)$ denotes the $\ell$th Fr\'echet differential of $f$ at $x$.

By \cite[Theorem~7.33.1]{Szego}, the Gegenbauer polynomials satisfy
$$|C_\ell^\lambda(a)|\leq C_\ell^\lambda(1)
=\frac{(2\lambda)_\ell}{\ell!},\qquad \lambda>0,\quad -1\leq a\leq1.$$
Taking $\lambda=\alpha/2$ gives
$$|\mathrm d^\ell f(x)[v,\ldots,v]|\leq(\alpha)_\ell R^{-\alpha-\ell}.$$

Since $\mathrm d^\ell f(x)$ is a continuous symmetric $\ell$-linear form on the real Hilbert space $\mathbb R^{n+1}$,
we can use the  Banach  theorem   \cite[Proposition~1.44]{Dineen} to obtain
$$\|\mathrm d^\ell f(x)\|=\sup_{|v|=1}|\mathrm d^\ell f(x)[v,\ldots,v]|,$$
where $\|\mathrm d^\ell f(x)\|$ denotes the operator norm of the $\ell$th Fr\'echet differential
$$\|\mathrm d^\ell f(x)\|=\sup_{\substack{v_1,\ldots,v_\ell\in\mathbb R^{n+1}\\ |v_1|=\cdots=|v_\ell|=1}}
\left|\mathrm d^\ell f(x)[v_1,\ldots,v_\ell]\right|.$$
Evaluating this multilinear form at the coordinate unit vectors, with the $j$th vector repeated $\beta_j$ times, gives
$$|\partial^\beta f(x)|\leq\|\mathrm d^\ell f(x)\|\leq(\alpha)_\ell|x|^{-\alpha-\ell}.$$

Let $\varepsilon_0=(1,0,\ldots,0)$. For $x=R\varepsilon_0$, $R>0$, and $\beta=(\ell,0,\ldots,0)$,
$$\left|\partial_{x_0}^{\ell}|x|^{-\alpha}\right|_{x=R\varepsilon_0}
=\left|\frac{d^\ell}{dR^\ell}R^{-\alpha}\right|
=(\alpha)_\ell R^{-\alpha-\ell},$$
which proves optimality of the uniform constant.

 Finally, from  the  estimate
\[\frac{(\alpha)_\ell}{\ell!}=\frac{\Gamma(\ell+\alpha)}{\Gamma(\alpha)\Gamma(\ell+1)}
\leq C_\alpha(1+\ell)^{\alpha-1},\qquad \ell\in\mathbb N,
\]
we can  prove  the last assertion.
\end{proof}
Since   $1\leq\frac{1+\ell}{\ell}\leq2,  \ell\geq1,$
we have
$$(1+\ell)^{\alpha-1}\leq 2^{\max\{\alpha-1,0\}}\ell^{\alpha-1}, \quad  \alpha>0,$$
which allows us to  recover and extend the derivative estimate in Proposition~\ref{c-d} as follows.
\begin{corollary}\label{prop:mixed-derivative-estimate}
For every $\alpha>0$, there exists $C_\alpha>0$, independent of $n$, $\beta$, and $x$, such that
\[|\partial^\beta |x|^{-\alpha}|\leq C_\alpha\ell^{\alpha-1}\ell!|x|^{-\alpha-\ell},
\qquad \ell=|\beta|\geq1,  x\in\mathbb R^{n+1}\setminus\{0\}.\]
\end{corollary}

\begin{theorem}\label{thm:kernel-series-convergence}
 Let $\lambda, s\in\mathbb R$ with $-1<\lambda<n$ and $-1<\lambda+2s<n$. Let  $\{a_\mu\}_{\mu\in\mathbb N^{n+1}}\subset\mathbb R_n$ satisfy
\[
r=\limsup_{k\to\infty}\Big(\sum_{|\mu|=k}\big(\frac{k!}{\mu!}\big)^2|a_\mu|^2\Big)^{1/(2k)}<+\infty.
\]
Then  the series
\[
\sum_{\mu\in\mathbb N^{n+1}}\bigl[(-\Delta)^sT_\mu^{(\lambda)}\bigr](x)a_\mu
\]
and all its derivative series converge absolutely and uniformly on compact subsets of $B(r,+\infty)=\{x\in\mathbb R^{n+1}:|x|>r\}$, with each term identified with its smooth representative away from the origin.
\end{theorem}
\begin{proof}Let $-1<\lambda<n$.
We first prove the assertion for $s=0$. Let $\beta\in\mathbb N^{n+1}$ and put $m=|\beta|$. Since $\overline{x}$ is a polynomial of degree one, the Leibniz rule gives
\[
\partial^\beta K_\lambda(x)
=\overline{x}\,\partial^\beta|x|^{-\lambda-2}
+\sum_{\substack{0\leq j\leq n\\ \beta_j>0}}
\beta_j\overline{e_j}\,\partial^{\beta-\varepsilon_j}|x|^{-\lambda-2},
\]
where  $\varepsilon_j=(0,\ldots,0,\underset{j\text{th}}{1},0,\ldots,0)\in\mathbb N^{n+1}$. Proposition~\ref{prop:sharp-mixed-derivative-estimate} therefore yields
\[
\begin{aligned}
|\partial^\beta K_\lambda(x)|
&\leq(\lambda+2)_m\big(1+\frac{m}{\lambda+m+1}\big)|x|^{-\lambda-1-m}\\
&\leq2(\lambda+2)_m|x|^{-\lambda-1-m}
\leq C_\lambda(1+m)^{\lambda+1}m!|x|^{-\lambda-1-m}.
\end{aligned}
\]
The estimate includes $m=0$, when the sum in the Leibniz formula is empty.

Fix $\nu\in\mathbb N^{n+1}$ and put $h=|\nu|$. For $|\mu|=k$, we obtain
\[
\begin{aligned}
|\partial^\nu T_\mu^{(\lambda)}(x)|
&\leq\frac{C_\lambda}{\mu!}(1+k+h)^{\lambda+1}(k+h)!|x|^{-\lambda-1-k-h}\\
&\leq C_{\lambda,h}(1+k)^{\lambda+1+h}\frac{k!}{\mu!}|x|^{-\lambda-1-k-h},
\end{aligned}
\]
where we use
\[
1+k+h\leq(1+h)(1+k),
\qquad (k+h)!=k!\prod_{j=1}^{h}(k+j)\leq(1+h)^h(1+k)^h k!,
\]
with the empty product interpreted as $1$ when $h=0$.

Set
$$A_k=\Big(\sum_{|\mu|=k}\big(\frac{k!}{\mu!}\big)^2|a_\mu|^2\Big)^{1/2}.$$
Let $E$ be a nonempty compact subset of $B(r,+\infty)$, put $R=\min_{x\in E}|x|$, and choose $r<\rho<R$. Since $\limsup_{k\to\infty}A_k^{1/k}=r$, there exists $C_\rho>0$ such that $A_k\leq C_\rho\rho^k$ for every $k\in\mathbb N$. The Cauchy--Schwarz inequality gives
\[
\sum_{|\mu|=k}\frac{k!}{\mu!}|a_\mu|
\leq\binom{k+n}{n}^{1/2}A_k
\leq(1+k)^{n/2}A_k.
\]
Each $\partial^\nu T_\mu^{(\lambda)}$ is paravector-valued. Hence Proposition~\ref{paravectors} and the preceding estimates imply
\[
\begin{aligned}
\sum_{|\mu|=k}\sup_{x\in E}|\partial^\nu T_\mu^{(\lambda)}(x)a_\mu|
&\leq C_{\lambda,h}(1+k)^{\lambda+1+h}R^{-\lambda-1-k-h}
\sum_{|\mu|=k}\frac{k!}{\mu!}|a_\mu|\\
&\leq C_{\lambda,h}C_\rho R^{-\lambda-1-h}
(1+k)^{\lambda+1+h+n/2}\left(\frac{\rho}{R}\right)^k.
\end{aligned}
\]
Since $\rho<R$, the right-hand side is summable over $k$, which proves absolute and uniform convergence on $E$ of the series and all its derivative series. This proves the assertion for $s=0$.

For general $s\in\mathbb R$ satisfying $-1<\lambda+2s<n$, Theorem~\ref{p-p+2s} gives
\[
(-\Delta)^sT_\mu^{(\lambda)}=c_{n,\lambda,s}T_\mu^{(\lambda+2s)}.
\]
Applying the case $s=0$ with $\lambda$ replaced by $\lambda+2s$ proves the assertion.
\end{proof}
We recall a version of Weierstrass theorem for partial-slice monogenic functions in \cite[Lemma 4.6]{Xu-Sabadini-4}.
\begin{lemma}\label{Weierstrass}
  Let $\{f_{m}\}_{m\in \mathbb{N}}$ be a sequence of  partial-slice monogenic functions of type $(p,q)$ in a partially symmetric slice
domain  $\Omega$  that converges to a function $f$ uniformly on compact sets in   $\Omega$. Then $f\in \mathcal{SM}_{(p,q)}(\Omega)$.
\end{lemma}

Now we are in a position to  establish a termwise Fueter--Sce--Qian theorem for partial-slice monogenic functions as follows.
\begin{theorem}\label{thm:termwise-fueter}
Let $p\in\{0,1,\ldots,n-1\}$, $q=n-p$, and  $s=j/2$, $j=0,1,\ldots,q-1$.  Suppose that $\{a_\mu\}_{\mu\in\mathbb N^{p+1}}\subset\mathbb R_n$ satisfies
$$r=\limsup_{k\to\infty}\Big(\sum_{|\mu|=k}\Big(\frac{k!}{\mu!}\Big)^2|a_\mu|^2\Big)^{1/(2k)}<+\infty.$$
Denote $f(x)=\sum_{\mu\in\mathbb N^{p+1}}T_\mu^{(p)}(x)a_\mu$ and define its termwise  fractional Laplacian transform
$$(-\Delta)^sf(x):=\sum_{\mu\in\mathbb N^{p+1}}\bigl[(-\Delta)^sT_\mu^{(p)}\bigr](x)a_\mu.$$
Then $f\in\mathcal{SM}_{(p,q)}(B(r,+\infty))$ and
\[(-\Delta)^sf\in\mathcal{SM}_{(p+2s,q-2s)}(B(r,+\infty)).
\]
\end{theorem}
\begin{proof}We identify each $\mu\in\mathbb N^{p+1}$ with $(\mu_0,\ldots,\mu_p,0,\ldots,0)\in\mathbb N^{n+1}$ and set all remaining coefficients $a_\mu$ equal to zero. This leaves the value of $r$ unchanged.
 For $s=j/2$, $j=0,1,\ldots,q-1$, we have $p+2s\in\{p,p+1,\ldots,n-1\}$.  Therefore, Theorems~\ref{p-p+2s} and~\ref{thm:kernel-series-convergence} give
\[
(-\Delta)^sf(x)=c_{n,p,s}\sum_{\mu\in\mathbb N^{p+1}}T_\mu^{(p+2s)}(x)a_\mu
\]
with convergence in $C^\infty_{\mathrm{loc}}(B(r,+\infty))$.

Note that $K_{p+2s}$ is partial-slice monogenic of type $(p+2s,q-2s)$ on $\mathbb R^{n+1}\setminus\{0\}$, so is every  finite sum
 $ \sum T_\mu^{(p+2s)}(x)a_\mu$ for $\mu\in\mathbb N^{p+1}$.
By Lemma \ref{Weierstrass},  we get
\[(-\Delta)^sf\in\mathcal{SM}_{(p+2s,q-2s)}(B(r,+\infty)),
\]
which completes the proof.
\end{proof}

\begin{remark}
Theorem~\ref{thm:termwise-fueter}, together with Theorem~\ref{thm:kernel-series-convergence}, extends the termwise fractional Laplacian construction for quaternionic Laurent series developed in \cite[Theorems~6.8 and~6.13]{Colombo-Martino} to series of mixed derivatives of Cauchy-type kernels in arbitrary dimension. In fact, for $p=0$, the coefficient condition becomes
$$r=\limsup_{k\to\infty}|a_k|^{1/k}<\infty.$$
Since the restriction of $T_{(k)}^{(0)}=T_{(k,0,\ldots,0)}^{(0)}=T_{-k-1}$ to $\mathbb R^{n+1}\setminus\{0\}$ is represented by $x^{-k-1}$, we obtain
$$f(x)=\sum_{k=0}^{\infty}x^{-k-1}a_k,\qquad |x|>r.$$
\end{remark}

\begin{remark}
For $s\in\mathbb N$, the termwise transform in Theorem~\ref{thm:termwise-fueter} agrees with the differential operator $(-\Delta)^s$ applied to $f$. For noninteger $s$, the theorem defines a termwise transform on $B(r,+\infty)$; it neither specifies a global tempered-distribution extension of $f$ nor identifies the resulting series with the Fourier-multiplier fractional Laplacian of such an extension.
\end{remark}

Finally, the following example shows that applying the fractional Laplacian to the zero extension need not preserve the Fueter mapping property.

\begin{example}\label{1-zero-extension}
Let $(p,q)=(0,2)$ and consider the slice monogenic function $f\equiv1$ on the unit ball $B_1\subset\mathbb R^3$. Let $F=\mathbf1_{B_1}$ be its zero extension. Since $F\in L^1(\mathbb R^3)$, its Fourier transform is bounded. Thus $2\pi|\xi|\widehat F(\xi)$ is locally integrable and has at most polynomial growth, so
\[
G:=(-\Delta)^{1/2}F=\mathcal F^{-1}\bigl[2\pi|\xi|\widehat F(\xi)\bigr]
\]
is well defined in $\mathscr S'(\mathbb R^3)$.

We claim that the restriction of $G$ to $B_1$ is represented by
\[
G(x)=\frac1{\pi^2}\int_{\mathbb R^3\setminus B_1}|x-y|^{-4}\,dy,\qquad x\in B_1.
\]
Indeed, for every $\varphi\in C_c^\infty(B_1)$, symmetry of the Fourier multiplier gives
\[
\langle G,\varphi\rangle=\int_{B_1}(-\Delta)^{1/2}\varphi(x)\,dx.
\]
The singular integral representation of the half-Laplacian (see e.g. \cite{Samko}) is
\[
(-\Delta)^{1/2}\varphi(x)=\frac1{\pi^2}\lim_{\varepsilon\downarrow0}\int_{|x-y|>\varepsilon}\frac{\varphi(x)-\varphi(y)}{|x-y|^4}\,dy.
\]
For each $\varepsilon>0$, antisymmetry under the interchange of $x$ and $y$ yields
\[
\int_{B_1}\int_{\substack{y\in B_1\\|x-y|>\varepsilon}}\frac{\varphi(x)-\varphi(y)}{|x-y|^4}\,dy\,dx=0.
\]
Moreover, $\varphi(y)=0$ outside $B_1$, and $\operatorname{supp}\varphi$ has positive distance from $\mathbb R^3\setminus B_1$. Consequently,
\[
\langle G,\varphi\rangle=\frac1{\pi^2}\int_{B_1}\varphi(x)\left(\int_{\mathbb R^3\setminus B_1}|x-y|^{-4}\,dy\right)dx.
\]
The passage to the limit is justified by the uniform convergence of the truncated singular integrals for $\varphi\in C_c^\infty(\mathbb R^3)$: cancellation of the linear Taylor term gives a truncation error of order $O(\varepsilon)$.

The displayed integral defining $G(x)$ is smooth on $B_1$, and differentiation under the integral is valid on every compact subset of $B_1$. Since $\Delta_x|x-y|^{-4}=12|x-y|^{-6}$, we obtain
\[
\Delta G(x)=\frac{12}{\pi^2}\int_{\mathbb R^3\setminus B_1}|x-y|^{-6}\,dy>0,\qquad x\in B_1.
\]
Hence $G$ is not harmonic and therefore not monogenic in $B_1$. This shows that applying the Fourier-multiplier fractional Laplacian to the zero extension does not generally yield the Fueter mapping property, even for a constant slice monogenic function.
\end{example}

\par
\noindent
\textbf{Acknowledgements}
 The   authors  would like to thank Professors   Fabrizio Colombo and Antonino De Martino for  their helpful  discussions.
\\
\textbf{Author contributions}
All authors have contributed equally to all aspects of this manuscript.
\\
\textbf{Conflict of interest}
There are no financial or non-financial interests that are directly or indirectly related to the work.
\\
\textbf{Data availability} 
Data sharing is not applicable to this article, as no datasets were generated during the study.

\vskip 10mm

\begingroup
\small
\begin{thebibliography}{99}
\setlength{\itemsep}{1pt}
\setlength{\parsep}{0pt}
\bibitem{Altavilla}
A. Altavilla,  H. De Bie, M. Wutzig, \textit{Implementing zonal harmonics with the Fueter principle}, J. Math. Anal. Appl. 495 (2021), no. 2, Paper No. 124764, 26 pp.

 \bibitem{Brackx} F. Brackx, R. Delanghe, F. Sommen, \textit{Clifford analysis}, Research Notes in Mathematics, Vol. 76, Pitman, Boston, 1982.

  \bibitem{Colombo-Sabadini-Sommen-10} F. Colombo, I. Sabadini, F. Sommen,     \textit{
The Fueter mapping theorem in integral form and the $\mathscr{F}$-functional calculus},
Math. Methods Appl. Sci. 33 (2010), no. 17, 2050-2066.

\bibitem{Colombo-Sabadini-Struppa-09} F. Colombo, I. Sabadini, D. C. Struppa, \textit{Slice monogenic functions}, Israel J. Math.  171   (2009), 385-403.

 \bibitem{Colombo-Martino} F. Colombo,  A. De Martino,   A. Debernardi Pinos,  \textit{Fractional $\frac{1}{2}$-Laplacian of holomorphic functions via Jacobi polynomials}, arXiv: 2607.22391, 2026.

     \bibitem{CDP25} F. Colombo, A. De Martino, S. Pinton,  \textit{Functions and operators of the polyharmonic and polyanalytic Clifford fine structures on the $S$-spectrum},  arXiv: 2501.14716, 2026.

\bibitem{Fivedim} F. Colombo, A. De Martino, S. Pinton, I. Sabadini, \textit{The fine structure of the spectral theory on the $S$-spectrum in dimension five}, J. Geom. Anal., 9, 33 (2023), 73 pp.


 \bibitem{Colombo-Martino-Qian-S} F. Colombo,  A. De Martino,  T. Qian, I. Sabadini,  \textit{The Poisson kernel and the Fourier transform of the slice monogenic Cauchy kernels}, J. Math. Anal. Appl. 512 (2022), no. 1, Paper No. 126115, 23 pp.






 \bibitem{De-Adan} A. De Martino,   A. Guzm\'{a}n Ad\'{a}n,   \textit{On the harmonic generalized Cauchy--Kovalevskaya extension and its connection with the Fueter--Sce theorem}, Trans. Amer. Math. Soc. 379 (2026), no. 8, 5391-5428.

\bibitem{Dineen} S. Dineen, \textit{Complex Analysis on Infinite Dimensional Spaces},
Springer Monographs in Mathematics, Springer-Verlag, London, 1999.

 \bibitem{Ding-Xu}    C. Ding, Z. Xu, \textit{Invariance of iterated global differential operator for slice monogenic
functions},   Comput. Methods Funct. Theory 25 (2025), no. 3, 735-752.
\bibitem{Dong-Qian-21} B. Dong,  T.  Qian,   \textit{Uniform generalizations of Fueter's theorem}, Ann. Mat. Pura Appl. (4) 200 (2021), no. 1, 229-251.
\bibitem{Fueter} R. Fueter,  \textit{Die funktionentheorie der differentialgleichungen $\Delta u = 0$ und $\Delta\Delta u = 0$ mit vier reellen variablen}. Commentarii Math. Helvetici 7 (1934), 307-330.
\bibitem{Eelbode} D. Eelbode,   \textit{The biaxial Fueter theorem},
Israel J. Math. 201 (2014), no. 1, 233-245.
\bibitem{Fei} M. Fei, P. Cerejeiras, U. K\"{a}hler,
\textit{Fueter's theorem and its generalizations in Dunkl-Clifford analysis},
J. Phys. A 42 (2009), no. 39, 395209, 15 pp.
\bibitem{Gentili-Struppa-07} G.  Gentili, D. C. Struppa, \textit{A new theory of regular functions of a quaternionic variable}, Adv. Math. 216 (2007) no. 1, 279-301.
\bibitem{Ghiloni-Stoppato-24-1} R. Ghiloni, C. Stoppato, \textit{A unified notion of regularity in one hypercomplex variable},  J. Geom. Phys. 202 (2024), Paper No. 105219, 13 pp.

\bibitem{Ghiloni-Stoppato-24-2} R. Ghiloni, C. Stoppato, \textit{A unified theory of regular functions of a hypercomplex variable},
Bull. Sci. Math.    209 (2026), Paper No. 103794, 92 pp.

\bibitem{Ghiloni-Stoppato-25} R. Ghiloni, C. Stoppato, \textit{A manifold Fueter--Sce phenomenon in one hypercomplex variable},
 J. Math. Anal. Appl. 564 (2026), no. 2, Paper No. 130903, 49 pp.
\bibitem{Gurlebeck} K. G\"{u}rlebeck, K. Habetha, W. Spr\"{o}{\ss}ig, \textit{Holomorphic functions in the plane and $n$-dimensional space}, Birkh\"{a}user Verlag, Basel, 2008.
\bibitem{Xie}S. He,  Y. Xie,  \textit{A new Fueter theorem for generalized partial-slice monogenic functions: direct and inverse results}, Mediterr. J. Math. 23 (2026), no. 6, Paper No. 181.
\bibitem{Hu-Ding} M. Hu, C. Ding, Y. Shen, J. Wang,   \textit{Integral formulas and Teodorescu transform for generalized partial-slice monogenic functions}, Adv. Appl. Clifford Algebr. 35 (2025), no. 5, Paper No. 50, 23 pp.

\bibitem{Huo-24} Q. Huo, P. Lian, J. Si, Z. Xu,  \textit{Almansi-type decomposition and Fueter--Sce theorem for
generalized partial-slice regular functions}, arXiv:2411.05571, 2024.

\bibitem{Huo-Sabadini-Xu-26}  Q. Huo, I. Sabadini,  Z. Xu, \textit{On the Fueter--Sce theorem and Cauchy--Kovalevskaya
extensions over alternative  $\ast$-algebras},  arXiv:2603.15014, 2026.

\bibitem{Kou} K. I. Kou,  T. Qian, F.  Sommen, \textit{Generalizations of Fueter's theorem},
Methods Appl. Anal. 9 (2002), no. 2, 273-289.

\bibitem{Perotti-22} A. Perotti, \textit{Cauchy--Riemann operators and local slice analysis over real alternative algebras}, J. Math. Anal. Appl. 516 (2022), no. 1, Paper No. 126480, 34 pp.
\bibitem{Perotti-26}   A. Perotti, \textit{Fueter trees for Dunkl-regular functions over alternative $\ast$-algebras},  arXiv:2604.12876, 2026.

\bibitem{Qian-97}  T. Qian, \textit{Generalization of Fueter's result to $\mathbb{R}^{n+1}$}. Atti Accad. Naz. Lincei Cl. Sci. Fis. Mat. Natur. Rend. Lincei (9) Mat. Appl. 8 (1997), no. 2, 111-117.
\bibitem{Qian-01}  T. Qian, \textit{Fourier analysis on starlike Lipschitz surfaces}, J. Funct. Anal. 183 (2001), no. 2, 370-412.
\bibitem{Qian-15}  T. Qian, \textit{Fueter Mapping Theorem in Hypercomplex Analysis}, Operator
Theory, Springer, Berlin (2015) pp. 1491-1507.

\bibitem{Rubin2025}
B. Rubin, \textit{Fractional Integrals, Potentials, and Radon Transforms}, 2nd ed., CRC Press, 2025.

\bibitem{Samko}
S. G. Samko, \textit{Hypersingular integrals and their applications. Analytical Methods and Special Functions,} 5. Taylor \& Francis Group, London, 2002.

\bibitem{Szego}
G. Szeg\H{o}, \textit{Orthogonal polynomials}, Fourth edition. American Mathematical Society Colloquium Publications, Vol. XXIII. American Mathematical Society, Providence, RI, 1975.

\bibitem{Sce} M. Sce,   \textit{Osservazioni sulle serie di potenze nei moduli quadratici},  Atti Accad. Naz. Lincei. Rend. Cl.
Sci. Fis. Mat. Nat. 23 (1957), 220-225.

\bibitem{Sommen-00} F. Sommen,   \textit{On a generalization of Fueter's theorem},
Z. Anal. Anwendungen 19 (2000), no. 4, 899-902.

\bibitem{Stein}E. Stein, \textit{Singular Integrals and Differentiability Properties of Functions}, Princeton University
Press, Princeton (1970).

 \bibitem{Xu-Sabadini} Z. Xu, I. Sabadini, \textit{Generalized partial-slice monogenic functions},  Trans. Amer. Math. Soc. 378 (2025), no. 2, 851-883.

\bibitem{Xu-Sabadini-3} Z. Xu, I. Sabadini, \textit{Generalized partial-slice monogenic functions: a synthesis of two function  theories}, Adv. Appl. Clifford Algebr. 34 (2024), no. 2, Paper No. 10.

 \bibitem{Xu-Sabadini-2} Z. Xu, I. Sabadini, \textit{On the Fueter--Sce theorem for generalized partial-slice monogenic functions},  Ann. Mat. Pura Appl.  204 (2025), no. 2, 835-857.

 \bibitem{Xu-Sabadini-4} Z. Xu, I. Sabadini, \textit{Segal--Bargmann transform for generalized partial-slice
monogenic functions},     Izv. Math. 89 (2025), no. 6, 1182-1207.

 \bibitem{Xu-Sabadini-26} Z. Xu, I. Sabadini, \textit{Generalized partial-slice monogenic functions: the octonionic case},  Trans. Amer. Math. Soc. 379 (2026), no. 5, 3257-3291.
\end{thebibliography}
\endgroup
\end{document}